\documentclass[12pt]{article}

\usepackage{amsmath, epsfig, cite, lineno}
\usepackage{amssymb,amsthm}
\usepackage{amsfonts, color}
\usepackage{latexsym}

\newtheorem{thm}{Theorem}[section]

\newtheorem{conj}[thm]{Conjecture}

\numberwithin{equation}{section}

\begin{document}


\begin{center}
{\Large\bf Factors of alternating sums of powers of\\[5pt] $q$-Narayana numbers}
\end{center}

\vskip 2mm \centerline{Victor J. W. Guo$^1$ and Qiang-Qiang Jiang$^2$}
\begin{center}
{\footnotesize $^1$School of Mathematical Sciences, Huaiyin Normal University, Huai'an, Jiangsu 223300,
 People's Republic of China\\
{\tt jwguo@hytc.edu.cn}\\[10pt]
$^2$Department of Mathematics, East China Normal University, Shanghai
200241,  People's Republic of China\\
{\tt 972126862@qq.com } }
\end{center}


\vskip 0.7cm \noindent{\bf Abstract.} The $q$-Narayana numbers $N_q(n,k)$ and $q$-Catalan numbers $C_n(q)$ are respectively defined  by
$$
N_q(n,k)=\frac{1-q}{1-q^n}{n\brack k}{n\brack k-1}\quad\text{and}\quad C_n(q)=\frac{1-q}{1-q^{n+1}}{2n\brack n},
$$
where ${n\brack k}=\prod_{i=1}^{k}\frac{1-q^{n-i+1}}{1-q^i}$.
We prove that, for any positive integers $n$ and $r$, there holds
\begin{align*}
\sum_{k=-n}^{n}(-1)^{k}q^{jk^2+{k\choose 2}}N_q(2n+1,n+k+1)^r \equiv 0 \pmod{C_n(q)},
\end{align*}
where $0\leqslant j\leqslant 2r-1$. We also propose several related conjectures.

\vskip 3mm \noindent {\it Keywords}: $q$-binomial coefficients, $q$-Narayana numbers, $q$-Catalan numbers

\vskip 0.2cm \noindent{\it AMS Subject Classifications:}  11A07, 05A30

\section{Introduction}

The Catalan numbers $C_n=\frac{1}{n+1}{2n\choose n}$ play an important role in combinatorics (see \cite{Stanley}).
It is well known that for any positive integer $n$,
$$
C_n=\sum_{k=1}^{n}N(n,k),
$$
where $N(n,k)=\frac{1}{n}{n\choose k}{n\choose k-1}$ are called Narayana numbers (see \cite{Narayana}).

In the past few years, many congruences on sums or alternating sums of binomial coefficients and
combinatorial numbers, such as Catalan numbers, Ap\'ery numbers, central Delannoy numbers, Schr\"oder numbers, Franel numbers,
have been obtained by Z.-W. Sun \cite{Sun1,Sun2,Sun3,Sun4,Sun5} and other authors \cite{GJZ,GZ,GZ2,Liu,Pan,PS,ST1,ST2}.

In this paper, motivated mainly by Z.-W. Sun's work, we shall prove the following
congruence on alternating sums of powers of Narayana numbers.

\begin{thm}\label{thm:main-1}
Let $n$ and $r$ be positive integers. Then
\begin{align}
\sum_{k=-n}^{n}(-1)^{k}N(2n+1,n+k+1)^r \equiv 0 \pmod{C_n}. \label{eq:main-1}
\end{align}
\end{thm}

We know that some congruences may have nice $q$-analogues (see, for example, \cite{GZ3,ShiPan,Pan0}). This is also the case for
the congruence \eqref{eq:main-1}. Recall that the $q$-shifted factorials (see \cite{Andrews98}) are defined by $(a;q)_0=1$
and $(a;q)_n=(1-a)(1-aq)\cdots (1-aq^{n-1})$ for $n=1,2,\ldots,$ and the $q$-binomial coefficients
are defined as
$$
{n\brack k}=
\begin{cases}\displaystyle\frac{(q;q)_n}{(q;q)_k (q;q)_{n-k}}, &\text{if $0\leqslant k\leqslant n$,} \\[10pt]
0, &\text{otherwise.}
\end{cases}
$$
For convenience, we let $[n]=\frac{1-q^n}{1-q}$ be a $q$-integer. It is natural to define the $q$-Narayana numbers $N_q(n,k)$ and the $q$-Catalan numbers $C_n(q)$ as follows:
$$
N_q(n,k)=\frac{1}{[n]}{n\brack k}{n\brack k-1},\quad C_n(q)=\frac{1}{[n+1]}{2n\brack n}.
$$
It is not difficult to see that both $q$-Narayana numbers and $q$-Catalan numbers are polynomials in $q$
with nonnegative integer coefficients (see \cite{Branden,FH}). Note that, the definition of $N_q(n,k)$ here differs by a
factor $q^{k(k-1)}$ from that in \cite{Branden}. We have the following $q$-analogue of Theorem \ref{thm:main-1}.

\begin{thm}\label{thm:main-2}
Let $n$ and $r$ be positive integers and let $0\leqslant j\leqslant 2r-1$. Then
\begin{align}
\sum_{k=-n}^{n}(-1)^{k}q^{jk^2+{k\choose 2}}N_q(2n+1,n+k+1)^r \equiv 0 \pmod{C_n(q)}.
\label{eq:main-2}
\end{align}
\end{thm}

It is easily seen that when $q\to 1$ the $q$-congruence \eqref{eq:main-2} reduces to \eqref{eq:main-1}.
It seems that \eqref{eq:main-2} also holds for $j\geqslant 2r$ (see \eqref{eq:fofk} in Conjecture \ref{conj:4} for a more general form).

\section{Proof of Theorem \ref{thm:main-2}}
Noticing that
\begin{align}
{2n+1\brack n+k}{2n+1\brack n+k+1}
=\frac{(q;q)_{2n+1}^2}{(q;q)_{2n}(q;q)_{2n+2}}{2n\brack n+k}{2n+2\brack n+k+1},
\label{eq:relation}
\end{align}
we can rewrite Theorem \ref{thm:main-2} in the following equivalent form.
\begin{thm}\label{thm:main-3}
Let $n$ and $r$ be positive integers and let $0\leqslant j\leqslant 2r-1$. Then
\begin{align}
&\hskip- 3mm \frac{(q;q)_{2n+1}^{2r}}{(q;q)_{2n}^{r}(q;q)_{2n+2}^{r}}
\sum_{k=-n}^{n}(-1)^{k}q^{jk^2+{k\choose 2}}{2n\brack n+k}^r{2n+2\brack n+k+1}^r \nonumber\\
&\equiv 0 \pmod{{2n+1\brack n}[2n+1]^{r-1} }.
\label{eq:main-3}
\end{align}
\end{thm}

In the paper \cite[Theorem 4.7]{GJZ}, Guo, Jouhet and Zeng proved the following result.
\begin{thm}\label{thm:GJZ}
For all positive integers $n_1,\ldots,n_m$ and $0\leqslant j\leqslant m-1$, the alternating sum
$$
(q;q)_{n_1}\prod_{i=1}^m\frac{(q;q)_{n_i+n_{i+1}}}{(q;q)_{2n_i}}
\sum_{k=-n_1}^{n_1}(-1)^k q^{jk^2+{k\choose 2}}\prod_{i=1}^m {2n_i\brack n_i+k},
$$
where $n_{m+1}=0$, is a polynomial in $q$ with nonnegative integer coefficients.
\end{thm}

In what follows, we will show that the congruence \eqref{eq:main-3} can be deduced from combining two special cases of
Theorem \ref{thm:GJZ}.

Denote the left-hand side of \eqref{eq:main-3} by $S(n,r,j)$. By the relation \eqref{eq:relation},
it is clear that $S(n,r,j)$ is a polynomial in $q$ with integer coefficients. Letting $m=2r$, $n_1=n_3=\cdots =n_{2r-1}=n$, and $n_2=n_4=\cdots=n_{2r}=n+1$ in
Theorem \ref{thm:GJZ}, we see that
\begin{align*}
\frac{(q;q)_{n}(q;q)_{n+1}(q;q)_{2n+1}^{2r-1}}{(q;q)_{2n}^{r}(q;q)_{2n+2}^{r}}
\sum_{k=-n}^{n}(-1)^{k}q^{jk^2+{k\choose 2}}{2n\brack n+k}^r{2n+2\brack n+k+1}^r
\in\mathbb{Z}[q],
\end{align*}
which can also be written as
\begin{align*}
{2n+1\brack n}^{-1}S(n,r,j)\in\mathbb{Z}[q].
\end{align*}
Namely,  $S(n,r,j)\equiv 0\pmod{{2n+1\brack n}}$, or
\begin{align}
[2n+1]^{r-1}S(n,r,j)\equiv 0\pmod{{2n+1\brack n}[2n+1]^{r-1}}. \label{eq:prime1}
\end{align}
On the other hand, letting $m=2r$, $n_1=n_2=\cdots =n_{r}=n$, and $n_{r+1}=n_{r+2}=\cdots=n_{2r}=n+1$ in
Theorem \ref{thm:GJZ}, we have
\begin{align*}
\frac{(q;q)_{n}(q;q)_{n+1}(q;q)_{2n+1}(q;q)_{2n}^{r-1}(q;q)_{2n+2}^{r-1}}{(q;q)_{2n}^{r}(q;q)_{2n+2}^{r}}
\sum_{k=-n}^{n}(-1)^{k}q^{jk^2+{k\choose 2}}{2n\brack n+k}^r{2n+2\brack n+k+1}^r
\in\mathbb{Z}[q],
\end{align*}
which,  by the relation $(q;q)_{2n}(q;q)_{2n+2}=(q;q)_{2n+1}^{2}\frac{[2n+2]}{[2n+1]}$, can be rewritten as
$$
{2n+1\brack n}^{-1}\frac{[2n+2]^{r-1}}{[2n+1]^{r-1}}S(n,r,j)\in\mathbb{Z}[q].
$$
Namely,
\begin{align}
[2n+2]^{r-1}S(n,r,j)\equiv 0\pmod{{2n+1\brack n}[2n+1]^{r-1}}. \label{eq:prime2}
\end{align}

It is easy to see that the polynomials $[2n+1]^{r-1}$ and $[2n+2]^{r-1}$ are relatively prime.
Therefore, by the Euclid algorithm for polynomials, there exist polynomials $P(q)$ and $Q(q)$ in $q$ with rational coefficients such that
\begin{align}
P(q)[2n+1]^{r-1}+Q(q)[2n+2]^{r-1}=1. \label{eq:prime3}
\end{align}
It follows from \eqref{eq:prime1}--\eqref{eq:prime3} that the congruence
\begin{align*}
S(n,r,j)\equiv 0\pmod{{2n+1\brack n}[2n+1]^{r-1}}.
\end{align*}
holds in the ring $\mathbb{Q}[q]$. In other words, there exists a polynomial $R(q)\in\mathbb{Q}[q]$
such that
\begin{align}
S(n,r,j)={2n+1\brack n}[2n+1]^{r-1}R(q).  \label{eq:final}
\end{align}
Since the polynomials $S(n,r,j)$ and ${2n+1\brack n}[2n+1]^{r-1}$ are in $\mathbb{Z}[q]$, and the leading coefficient of the latter is one, the identity \eqref{eq:final} means that $R(q)\in\mathbb{Z}[q]$. This completes the proof.

\section{Some open problems}
It seems that Theorems \ref{thm:main-1} and \ref{thm:main-2} can be further generalized as follows.
\begin{conj}\label{conj:1}
Let $n_1,\ldots,n_m$ be positive integers. Then
\begin{align*}
&\hskip -3mm \sum_{k=-n_1}^{n_1}(-1)^k \prod_{i=1}^m {n_i+n_{i+1}+1\choose n_i+k}{n_i+n_{i+1}+1\choose n_i+k+1} \\
&\equiv 0 \pmod{{n_1+n_m+1\choose n_1}\prod_{i=1}^{m-1}(n_i+n_{i+1}+1)},
\end{align*}
where $n_{m+1}=n_1$.
\end{conj}

\begin{conj}\label{conj:2}
Let $n$ and $r$ be positive integers and let $0\leqslant j\leqslant 2r-1$. Then
\begin{align*}
\frac{1}{C_n(q)}\sum_{k=-n}^{n}(-1)^{k}q^{jk^2+{k\choose 2}}N_q(2n+1,n+k+1)^r
\end{align*}
is a polynomial in $q$ with nonnegative integer coefficients.
\end{conj}

Note that the upper bound $2r-1$ of $j$ in Conjecture \ref{conj:2} seems to be the best possible.
Numerical calculation implies that Conjecture \ref{conj:2} does not hold when $j\geqslant 2r$.
Furthermore, we have the following generalization of Conjecture \ref{conj:2}.

\begin{conj}\label{conj:3}
For all positive integers $n_1,\ldots,n_m$ and $0\leqslant j\leqslant 2m-1$, the expression
$$
{n_1+n_m+1\brack n_1}^{-1}\prod_{i=1}^{m-1}\frac{1}{[n_i+n_{i+1}+1]}\sum_{k=-n_1}^{n_1}(-1)^k q^{jk^2+{k\choose 2}}
\prod_{i=1}^m {n_i+n_{i+1}+1\brack n_i+k}{n_i+n_{i+1}+1\brack n_i+k+1},
$$
where $n_{m+1}=n_1$, is a polynomial in $q$ with nonnegative integer coefficients.
\end{conj}

We end the paper with the following $q$-analogue of Conjecture \ref{conj:1}.
\begin{conj}\label{conj:4}
Let $n_1,\ldots,n_m$ be positive integers, and let $f(k)$ be a polynomial in $k$ with
integer coefficients. Then
\begin{align*}
&\hskip -3mm \sum_{k=-n_1}^{n_1}(-1)^k q^{f(k)+{k\choose 2}}\prod_{i=1}^m {n_i+n_{i+1}+1\brack n_i+k}{n_i+n_{i+1}+1\brack n_i+k+1} \\
&\equiv 0 \pmod{{n_1+n_m+1\brack n_1}\prod_{i=1}^{m-1}[n_i+n_{i+1}+1]},
\end{align*}
where $n_{m+1}=n_1$. In particular, for any positive integer $n$, we have
\begin{align}
\sum_{k=-n}^{n}(-1)^{k}q^{f(k)+{k\choose 2}}N_q(2n+1,n+k+1)^r \equiv 0 \pmod{C_n(q)}. \label{eq:fofk}
\end{align}
\end{conj}

\vskip 5mm \noindent{\bf Acknowledgments.} This work was partially
supported by the National Natural Science Foundation of China (grant 11371144),
the Natural Science Foundation of Jiangsu Province (grant BK20161304),
and the Qing Lan Project of Education Committee of Jiangsu Province.


\begin{thebibliography}{99}
\small \setlength{\itemsep}{-.8mm}

\bibitem{Andrews98}G.E. Andrews, The Theory of Partitions, Cambridge University Press, Cambridge, 1998.

\bibitem{Branden}P. Br\"and\'en, $q$-Narayana numbers and the flag $h$-vector of $J(2\times n)$, Discrete Math.
281 (2004), 67--81.

\bibitem{FH}J. F\"urlinger and J. Hofbauer, $q$-Catalan numbers, J. Combin. Theory, Ser. A 2 (1985), 248--264.

\bibitem{GJZ}V.J.W. Guo, F. Jouhet, J. Zeng, Factors of alternating sums
of products of binomial and $q$-binomial coefficients, Acta Arith. 127 (2007), 17--31.

\bibitem{GZ}V.J.W. Guo and J. Zeng, New congruences for sums involving Ap\'ery numbers or central Delannoy numbers,
preprint, Int. J. Number Theory 8 (2012), 2003--2016.

\bibitem{GZ2}V.J.W. Guo and J. Zeng, Proof of some conjectures of Z.-W. Sun on congruences for Ap\'ery polynomials,
J. Number Theory 132 (2012), 1731--1740.

\bibitem{GZ3}V.J.W. Guo and J. Zeng, Some $q$-analogues of supercongruences of Rodriguez-Villegas,
J. Number Theory 145 (2014), 301--316.

\bibitem{Liu}J.-C. Liu, A supercongruence involving Delannoy numbers and Schr\"oder numbers,
J. Number Theory  168  (2016), 117--127.

\bibitem{Narayana}T.V. Narayana, A partial order and its application to probability,
Sankhy\={a} 21 (1959), 91--98.


\bibitem{Pan0}H. Pan, A $q$-analogue of Lehmer's congruence, Acta Arith. 128 (2007), 303--318.

\bibitem{Pan}H. Pan, On divisibility of sums of Ap\'ery polynomials, J. Number Theory 143 (2014), 214--223.

\bibitem{PS}H. Pan and Z.-W. Sun, A combinatorial identity with application to Catalan numbers,
Discrete Math. 306 (2006), 1921--1940.

\bibitem{ShiPan}L.-L. Shi and H. Pan, A $q$-analogue of Wolstenholme's harmonic series congruence,
 Amer. Math. Monthly 114 (2007), 529¨C531.

\bibitem{Stanley}R.P. Stanley, Enumerative Combinatorics, Vol.~II, Cambridge University Press, Cambridge, 1999.

\bibitem{Sun1}Z.-W. Sun,  On Delannoy numbers and Schr\"oder numbers, J. Number Theory 131 (2011), 2387--2397.

\bibitem{Sun1}Z.-W. Sun, On sums of Ap\'ery polynomials and related congruences, J. Number Theory, 132 (2012),
2673--2699.

\bibitem{Sun2}Z.-W. Sun, Connections between $p=x^2+3y^2$ and Franel numbers, J. Number Theory 133 (2013), 2914--2928.

\bibitem{Sun3}Z.-W. Sun, Supercongruences involving products of two binomial coefficients,
¡¡¡¡Finite Fields Appl. 22 (2013), 24-44.

\bibitem{Sun4}Z.-W. Sun, Congruences for Franel numbers, Adv. Appl. Math. 51 (2013),  524--535.

\bibitem{Sun5}Z.-W. Sun, Congruences involving generalized trinomial coefficients, Sci. China 57 (2014), 1375--1400.

\bibitem{ST1}Z.-W. Sun and R. Tauraso, New congruences for central binomial coefficients,
Adv. Appl. Math. 45 (2010), 125--148.

\bibitem{ST2}Z.-W. Sun and R. Tauraso, On some new congruences for binomial coefficients,
¡¡¡¡Int. J. Number Theory 7(2011), 645--662.


\end{thebibliography}
\end{document}